\documentclass[11pt]{article}
\usepackage[margin=1in]{geometry}
\usepackage{amsmath,amssymb,amsthm,mathtools}
\usepackage{bm}
\usepackage{microtype}
\usepackage{authblk}
\usepackage[hidelinks]{hyperref}

\newtheorem{theorem}{Theorem}[section]
\newtheorem{lemma}[theorem]{Lemma}

\newcommand{\R}{\mathbb{R}}
\newcommand{\Sn}{\mathbb{S}^{n}}
\newcommand{\Ga}{\Gamma}
\newcommand{\tGa}{\widetilde\Gamma}
\newcommand{\sig}{\sigma}
\newcommand{\al}{\alpha}
\newcommand{\ft}{\widetilde f}

\title{Curvature Estimates and Prescribed Curvature Problems for Sum Hessian Equations on Star-Shaped Hypersurfaces}
\author[1]{Dongjing Ma\thanks{E-mail: \href{mailto:3261019581@qq.com}{3261019581@qq.com}}}
\author[1]{Peihe Wang\thanks{E-mail: \href{mailto:peihewang@hotmail.com}{peihewang@hotmail.com}}}
\author[2]{Hua Zhu\thanks{E-mail: \href{mailto:zhuhmaths@mail.ustc.edu.cn}{zhuhmaths@mail.ustc.edu.cn}}}

\affil[1]{School of Mathematical Sciences, Qufu Normal University, Qufu 273165, Shandong Province, People's Republic of China}
\affil[2]{School of Mathematics and Physics, Southwest University of Science and Technology, Mianyang 621010, Sichuan Province, People's Republic of China}

\date{}

\begin{document}
\maketitle

\begin{abstract}
Let $M$ be a closed star-shaped hypersurface in $\mathbb{R}^{n+1}$ and let
$\eta=Hg-h$ be the first Newton transformation of its second fundamental form. For a fixed
constant $\alpha>0$ and $2\le k\le n-1$, we study the prescribed sum Hessian curvature equation
\[
    \sigma_k(\lambda(\eta))+\alpha\sigma_{k-1}(\lambda(\eta))=f(X,\nu).
\]
Under the natural admissibility condition
\[
    \lambda(\eta)\in\Gamma_{k-1}\cap
    \{\lambda:\sigma_k(\lambda)+\alpha\sigma_{k-1}(\lambda)>0\},
\]
we establish a global curvature estimate for closed star-shaped solutions with a general positive
right-hand side depending on the position vector and the unit normal. A key point is that the sum
Hessian operator is nonhomogeneous; the resulting lower-order term in the maximum-principle
calculation is controlled by a partial curvature estimate. We also use the augmented-variable
identity
\[
    \sigma_k(\lambda)+\alpha\sigma_{k-1}(\lambda)
    =\sigma_k(\lambda_1,\ldots,\lambda_n,\alpha)
\]
to recover the relevant ellipticity and concavity properties. As an application of the a priori
estimates, we prove existence and uniqueness for the corresponding prescribed curvature problem
under suitable spherical barrier and radial monotonicity assumptions.
\end{abstract}



\section{Introduction}

Let $M$ be a hypersurface in the $(n+1)$-dimensional Euclidean space $\R^{n+1}$. Let
$\kappa(X)$, $H(X)$, and $\nu(X)$ denote the principal curvatures, mean curvature, and unit
outward normal vector of $M$ at $X$, respectively. Define the $(0,2)$-tensor
\[
    \eta_{ij}=Hg_{ij}-h_{ij},
\]
where $g_{ij}$ and $h_{ij}$ are the coefficients of the first and second fundamental forms.
The tensor $\eta$ is the first Newton transformation of $h$ with respect to $g$. If
$\lambda(\eta)=(\lambda_1,\ldots,\lambda_n)$ denotes the eigenvalues of $\eta$, then
\[
    \lambda_i=H-\kappa_i=\sum_{j\ne i}\kappa_j,
    \qquad i=1,\ldots,n.
\]

The classical prescribed curvature equation
\[
    \sig_k(\kappa(X))=f(X,\nu(X))
\]
plays a fundamental role in fully nonlinear elliptic equations and geometric analysis. The cases
$k=1$, $k=2$, and $k=n$ correspond, up to the usual normalization, to the prescribed mean,
scalar, and Gauss curvature equations. The classical theory of fully nonlinear Hessian equations and prescribed curvature problems on star-shaped hypersurfaces can be found in \cite{CNS1985,CNS1986}. In a different but closely related direction, Chu and
Jiao~\cite{ChuJiao2021} studied the Hessian-type curvature equation
\[
    \sig_k(\lambda(\eta))=f(X,\nu(X))
\]
for star-shaped hypersurfaces. They established a global curvature estimate under the natural
condition $\lambda(\eta)\in\Ga_k$ and applied it to the corresponding prescribed curvature
problem. Their proof combines a maximum-principle test function with a partial curvature
estimate that provides the good third-order terms needed to control the derivatives of the
right-hand side. Related global curvature estimates for convex hypersurfaces and prescribed curvature measure problems were established in \cite{GRW2015,GLL2012}. The evolution of star-shaped hypersurfaces by symmetric functions of their principal curvatures was also studied in \cite{Urbas1990}.

Sum Hessian equations constitute another natural class of fully nonlinear equations. In Euclidean
domains the model equation is
\[
    \sig_k(D^2u)+\al\sig_{k-1}(D^2u)=f(x,u,Du),
    \qquad \al>0.
\]
The operator is a linear combination of two neighboring Hessian operators of different degrees
and is therefore nonhomogeneous. Li, Ren, and Wang~\cite{LiRenWang2019} identified the natural
admissible region
\[
    \Ga_{k-1}\cap\{\lambda:\sig_k(\lambda)+\al\sig_{k-1}(\lambda)>0\}
\]
and studied curvature estimates for related linear combinations of elementary symmetric
functions on convex hypersurfaces. Liu and Ren~\cite{LiuRen2023} obtained Pogorelov-type
$C^2$ estimates for the standard sum Hessian equation. Liang, Yan, and Zhu~\cite{LiangYanZhu2026}
subsequently developed Pogorelov estimates and rigidity results for $(k-1)$-convex
solutions under an additional semi-convexity assumption, stressing that the loss of homogeneity prevents a direct transfer of
several arguments from the pure $k$-Hessian equation. Further developments include the
Pogorelov estimates of Li and Ren~\cite{LiRen2025} and the interior estimate of Ren and Sun~\cite{RenSun2026} for a two-sum
Hessian equation with a nonnegative variable coefficient $\alpha(x)\ge 0$ satisfying the
structural condition $D^2\alpha=\beta I$, under the additional assumption $\Delta u>0$.

Of particular relevance to the present problem is the work of Ren and Wang~\cite{RenWang2025},
who considered
\[
    \sig_k(\mu)+\al\sig_{k-1}(\mu)=f(x,u,Du),
    \qquad
    \mu_i=\sum_{j\ne i}\lambda_j(D^2u),
\]
so that $\mu$ is the spectrum of $(\Delta u)I-D^2u$. They established interior $C^2$ and
Pogorelov-type estimates for $0<k<n$, while the endpoint $k=n$ requires a weaker estimate.
This is algebraically the closest Euclidean-domain analogue of our operator. The geometric
setting here is different: $\eta=Hg-h$ is defined by the second fundamental form of a closed
hypersurface, and the goal is a global curvature estimate rather than an interior Hessian
estimate.

There are also recent geometric works on sum curvature equations formulated directly in the
principal curvatures. For example, Yang and Lu~\cite{YangLu2025} study complete hypersurfaces
of constant sum Hessian curvature in hyperbolic space, with curvature function
$\sig_{n-1}(\kappa)+\al\sig_n(\kappa)$. These problems, as well as the convex-hypersurface
results of Li, Ren, and Wang~\cite{LiRenWang2019}, are structurally distinct from applying the
sum Hessian operator to the eigenvalues of the first Newton transformation $Hg-h$.

Motivated by the above developments and by the geometric framework of Chu and Jiao, we study
the prescribed sum Hessian curvature equation associated with $\eta$. Fix a constant $\al>0$.
For $1\le j\le n$, set
\[
    S_{j,\al}(\lambda)=\sig_j(\lambda)+\al\sig_{j-1}(\lambda),
    \qquad \sig_0\equiv1.
\]
Our equation is
\begin{equation}\label{eq:main}
    S_{k,\al}(\lambda(\eta))
    =\sig_k(\lambda(\eta))+\al\sig_{k-1}(\lambda(\eta))
    =f(X,\nu(X)),\qquad X\in M.
\end{equation}
The natural admissible region is
\[
    \tGa_{k,\al}
    :=\Ga_{k-1}\cap\{\lambda\in\R^n:S_{k,\al}(\lambda)>0\},
\]
where
\[
    \Ga_j=\{\lambda\in\R^n:\sig_1(\lambda)>0,\ldots,\sig_j(\lambda)>0\}.
\]
A $C^2$ hypersurface is called $(\eta,k,\al)$-admissible if
$\lambda(\eta)\in\tGa_{k,\al}$ everywhere on $M$.

The key structural observation is obtained by introducing the augmented vector
\[
    \widehat\lambda=(\lambda_1,\ldots,\lambda_n,\al)\in\R^{n+1}.
\]
Then
\[
    S_{j,\al}(\lambda)=\sig_j(\widehat\lambda).
\]
Thus the nonhomogeneous sum Hessian operator is the restriction of the standard Hessian operator
in one higher dimension to the affine hyperplane whose last coordinate equals $\al$. This
recovers ellipticity and concavity from $\sig_k^{1/k}$ in one higher dimension. On the other
hand, the loss of homogeneity changes the Euler identity and produces an additional lower-order
term in the maximum-principle calculation. The main analytical task is to absorb this term using
the partial curvature estimate.

To the best of our knowledge, a global curvature estimate and a corresponding prescribed-curvature
existence theorem for \eqref{eq:main} on closed star-shaped hypersurfaces under the natural
$(\eta,k,\al)$-admissibility condition have not previously been established. Our first main
result is the following.

\begin{theorem}\label{thm:curvature}
Let $2\le k\le n-1$ and let $\al>0$. Let $M$ be a closed, star-shaped,
$(\eta,k,\al)$-admissible hypersurface in $\R^{n+1}$ satisfying \eqref{eq:main}, where
$f>0$ and $f\in C^2(\mathcal U)$ for an open neighborhood $\mathcal U$ of the unit normal
bundle of $M$ in $\R^{n+1}\times\Sn$. Then
\[
    \max_{X\in M}\max_{1\le i\le n}|\kappa_i(X)|\le C,
\]
where
\[
    C=C\bigl(n,k,\al,\|M\|_{C^1},\inf f,\|f\|_{C^2}\bigr).
\]
\end{theorem}

The restriction $k\le n-1$ enters in the partial curvature estimate in
Section~\ref{sec:curvature}; the endpoint $k=n$ requires a different argument and is not treated
here. The proof follows the maximum-principle strategy of Chu and Jiao, but the lower-order
$\al\sig_{k-1}$ term changes the Euler identity and gives rise to an additional quantity
$R_\al$. The principal new point in the estimate is to control this term without strengthening
the natural admissibility condition.

We also prove an existence result. Put
\begin{equation}\label{eq:phi-alpha}
    \Phi_\al(r)
    :=\binom nk\frac{(n-1)^k}{r^k}
      +\al\binom n{k-1}\frac{(n-1)^{k-1}}{r^{k-1}}.
\end{equation}
This is the value of the left-hand side of \eqref{eq:main} on the sphere of radius $r$. Since the
sum Hessian operator is not homogeneous, the spherical barriers used for the pure $k$-Hessian
case acquire the second term in \eqref{eq:phi-alpha}.

\begin{theorem}\label{thm:existence}
Let $2\le k\le n-1$, $\al>0$, and $0<r_1<1<r_2$. Let
\[
  f\in C^2\bigl((\overline B_{r_2}\setminus B_{r_1})\times\Sn\bigr)
\]
be positive and assume
\begin{equation}\label{eq:barrier}
\begin{aligned}
  f\left(X,\frac{X}{|X|}\right)&\ge \Phi_\al(r_1),&& |X|=r_1,\\
  f\left(X,\frac{X}{|X|}\right)&\le \Phi_\al(r_2),&& |X|=r_2.
\end{aligned}
\end{equation}
Assume moreover that for every fixed $x,\nu\in\Sn$,
\begin{equation}\label{eq:radial-monotonicity}
    \frac{\partial}{\partial\rho}\bigl[\rho^k f(\rho x,\nu)\bigr]\le0,
    \qquad r_1\le\rho\le r_2.
\end{equation}
Then, for every $0<\beta<1$, equation \eqref{eq:main} admits a unique
$C^{3,\beta}$ star-shaped $(\eta,k,\al)$-admissible solution satisfying
$r_1\le |X|\le r_2$.
\end{theorem}

When $\al=0$, equation \eqref{eq:main} formally reduces to the Hessian-type curvature equation
of Chu and Jiao. For $\al>0$, both the global curvature estimate and the existence argument must
account for the nonhomogeneous lower-order Hessian term. In particular, the curvature estimate
requires a separate control of $R_\al$, while the continuity method uses the modified spherical
value \eqref{eq:phi-alpha}.

The paper is organized as follows. Section~\ref{sec:preliminaries} records the structural
properties of the sum Hessian operator. Section~\ref{sec:curvature} proves the global curvature
estimate. Section~\ref{sec:existence} establishes the a priori estimates and completes the
existence and uniqueness proof by the continuity method.

\section{Preliminaries}\label{sec:preliminaries}

We first recall some basic definitions and properties of the elementary symmetric functions.
For $1\le k\le n$, the $k$-th elementary symmetric function is defined by
\[
    \sig_k(\lambda)
    =
    \sum_{1\le i_1<\cdots<i_k\le n}
    \lambda_{i_1}\cdots\lambda_{i_k},
\]
for $\lambda=(\lambda_1,\lambda_2,\ldots,\lambda_n)\in\R^n$.
For convenience, we further set
\[
    \sig_0(\lambda)=1,
    \qquad
    \sig_k(\lambda)=0
    \quad\text{for }k>n\text{ or }k<0.
\]

We next record the structural properties of the Hessian sum operator.

\begin{lemma}\label{lem:augmented}
Let $\al>0$ and let $\widehat\lambda=(\lambda,\al)\in\R^{n+1}$. Then
\[
    S_{j,\al}(\lambda)=\sig_j(\widehat\lambda),\qquad 1\le j\le n.
\]
Moreover,
\[
    \lambda\in\tGa_{k,\al}
    \quad\Longleftrightarrow\quad
    \widehat\lambda\in\Ga_k^{n+1}.
\]
Consequently, $\tGa_{k,\al}$ is open and convex, and
\[
    G(\lambda):=S_{k,\al}(\lambda)^{1/k}
\]
is elliptic and concave on $\tGa_{k,\al}$.
\end{lemma}

\begin{proof}
The identity
\[
  \sig_j(\lambda_1,\ldots,\lambda_n,\al)
  =\sig_j(\lambda)+\al\sig_{j-1}(\lambda)
\]
is immediate from the definition of the elementary symmetric functions. If
$\lambda\in\tGa_{k,\al}$, then for $1\le j\le k-1$,
\[
  \sig_j(\widehat\lambda)=\sig_j(\lambda)+\al\sig_{j-1}(\lambda)>0,
\]
and $\sig_k(\widehat\lambda)=S_{k,\al}(\lambda)>0$. Hence
$\widehat\lambda\in\Ga_k^{n+1}$.

Conversely, if $\widehat\lambda\in\Ga_k^{n+1}$, then for
$1\le j\le k-1$ the standard ellipticity of $\sig_{j+1}$ on $\Ga_{j+1}$ gives
\[
   \sig_j(\lambda)
   =\frac{\partial\sig_{j+1}(\widehat\lambda)}{\partial\widehat\lambda_{n+1}}>0.
\]
Thus $\lambda\in\Ga_{k-1}$, while
$S_{k,\al}(\lambda)=\sig_k(\widehat\lambda)>0$.
The convexity statement follows by intersecting the convex cone $\Ga_k^{n+1}$ with the affine
hyperplane $\widehat\lambda_{n+1}=\al$. Finally,
$G(\lambda)=\sig_k(\widehat\lambda)^{1/k}$, so ellipticity and concavity follow from the
corresponding standard properties of $\sig_k^{1/k}$ on $\Ga_k$~\cite{CNS1985}.
\end{proof}

For later use, if $\lambda$ is diagonal, write
\[
  S_k^{ii}=\frac{\partial S_{k,\al}}{\partial\lambda_i},
  \qquad
  G^{ii}=\frac{\partial G}{\partial\lambda_i}.
\]
Then
\begin{equation*}
  S_k^{ii}=\sig_{k-1}(\lambda|i)+\al\sig_{k-2}(\lambda|i),
\end{equation*}
\begin{equation}\label{eq:trace-S}
  \sum_iS_k^{ii}
  =(n-k+1)\sig_{k-1}(\lambda)+\al(n-k+2)\sig_{k-2}(\lambda),
\end{equation}
and
\begin{equation*}
  \sum_i\lambda_iS_k^{ii}
  =kS_{k,\al}(\lambda)-\al\sig_{k-1}(\lambda).
\end{equation*}
Thus, on a solution of \eqref{eq:main}, if
\begin{equation}\label{eq:Ralpha}
  \ft=f^{1/k},
  \qquad
  R_\al:=\frac{\al}{k}f^{1/k-1}\sig_{k-1}(\lambda),
\end{equation}
then
\begin{equation*}
  \sum_iG^{ii}\lambda_i=\ft-R_\al.
\end{equation*}

For $F(\kappa)=G(\eta(\kappa))$, at a fixed point, choose a local orthonormal frame such that $g_{ij}=\delta_{ij}$, $h_{ij}=\kappa_i\delta_{ij}$, then
\begin{equation}\label{eq:Fii}
  F^{ii}=\frac{\partial F}{\partial\kappa_i}
  =\sum_{p\ne i}G^{pp}.
\end{equation}
Consequently,
\begin{equation}\label{eq:trace-Fh}
  \sum_iF^{ii}h_{ii}=\sum_iG^{ii}\eta_{ii}=\ft-R_\al.
\end{equation}

\begin{lemma}\label{lem:trace-lower}
For $\lambda\in\tGa_{k,\al}$,
\begin{equation}\label{eq:trace-lower}
    \sum_iF^{ii}\ge c_{n,k}>0,
\end{equation}
and
\begin{equation}\label{eq:Ralpha-trace}
    0\le R_\al\le
    \frac{\al}{(n-1)(n-k+1)}\sum_iF^{ii}.
\end{equation}
\end{lemma}

\begin{proof}
By \eqref{eq:Fii} and \eqref{eq:trace-S},
\begin{align*}
 \sum_iF^{ii}
 &= (n-1)\sum_iG^{ii}\\
 &=\frac{n-1}{k}f^{1/k-1}
 \bigl[(n-k+1)\sig_{k-1}(\lambda)
 +\al(n-k+2)\sig_{k-2}(\lambda)\bigr]\\
 &\ge \frac{(n-1)(n-k+1)}{k}f^{1/k-1}
 \sig_{k-1}(\widehat\lambda).
\end{align*}
Since $\widehat\lambda\in\Ga_k^{n+1}$, the Newton--Maclaurin inequality gives
\[
  \sig_{k-1}(\widehat\lambda)\ge c_{n,k}f^{(k-1)/k},
\]
which proves \eqref{eq:trace-lower}. Formula \eqref{eq:Ralpha-trace} follows immediately by
comparing \eqref{eq:Ralpha} with the displayed formula for $\sum_iF^{ii}$.
\end{proof}

We now recall the standard geometric formulas for star-shaped radial graphs that will be used in the sequel.

Let $M$ be the radial graph of a positive function $\rho$ on $\Sn$,
\[
    X(x)=\rho(x)x,\qquad x\in\Sn.
\]
If $\widehat g$ and $\nabla$ denote the standard metric and covariant derivative on $\Sn$, then
\begin{equation}\label{eq:radial-formulas}
    g_{ij}=\rho^2\widehat g_{ij}+\rho_i\rho_j,
    \qquad
    h_{ij}=\frac{\rho^2\widehat g_{ij}+2\rho_i\rho_j-\rho\rho_{ij}}
    {\sqrt{\rho^2+|\nabla\rho|^2}},
\end{equation}
and
\begin{equation*}
    \nu=\frac{\rho x-\nabla\rho}{\sqrt{\rho^2+|\nabla\rho|^2}}.
\end{equation*}
In particular, with
\begin{equation*}
    u=\langle X,\nu\rangle,
\end{equation*}
we have
\begin{equation}\label{eq:support-rho}
    u=\frac{\rho^2}{\sqrt{\rho^2+|\nabla\rho|^2}}.
\end{equation}

Fix $X_0\in M$ and let $\{e_1,\ldots,e_n\}$ be a local orthonormal frame near $X_0$.
With the sign convention used throughout the paper,
\begin{align}
    X_{ij}&=-h_{ij}\nu,\label{eq:gauss}\\
    \nu_i&=h_{ij}e_j,\notag\\
    h_{ijk}&=h_{ikj},\label{eq:codazzi}\\
    R_{ijkl}&=h_{ik}h_{jl}-h_{il}h_{jk}.\notag
\end{align}
We shall also use
\begin{equation}\label{eq:commutation-h}
 h_{ijkl}=h_{klij}+(h_{mj}h_{il}-h_{ml}h_{ij})h_{mk}
 +(h_{mj}h_{kl}-h_{ml}h_{kj})h_{mi}.
\end{equation}

\section{Curvature estimate}\label{sec:curvature}

We prove Theorem~\ref{thm:curvature}. Since $\lambda(\eta)\in\Ga_{k-1}\subset\Ga_1$,
\[
   0<\sig_1(\lambda(\eta))=(n-1)H,
\]
so $H>0$. It is therefore enough to obtain a uniform upper bound for the largest principal
curvature $\kappa_{\max}$.

By \eqref{eq:support-rho} there is a constant $C$, depending only on the $C^1$ geometry of the
radial graph, such that
\begin{equation}\label{eq:u-bounds}
   C^{-1}\le u\le C.
\end{equation}
Set
\[
    a=\frac12\inf_M u>0.
\]
On the set $D=\{\kappa_{\max}>0\}$ consider
\begin{equation*}
   Q=\log\kappa_{\max}-\log(u-a)+\frac A2|X|^2,
\end{equation*}
where $A>1$ will be fixed later. The function $Q$ tends to $-\infty$ at $\partial D$ and hence
attains its maximum at some $X_0\in D$. Choose a local orthonormal frame near $X_0$ such that
\[
    h_{ij}=\delta_{ij}h_{ii},
    \qquad
    h_{11}\ge h_{22}\ge\cdots\ge h_{nn}
\]
at $X_0$. Then
\[
    \eta_{11}\le\eta_{22}\le\cdots\le\eta_{nn}.
\]
If $\kappa_{\max}(X_0)$ has multiplicity, we use the standard Calabi perturbation of the
second fundamental form so that the largest eigenvalue becomes simple at $X_0$.  The perturbation
may be chosen arbitrarily small; all estimates below are uniform in the perturbation parameter, and
we let that parameter tend to zero at the end.  Thus, without loss of generality for the maximum-
principle calculation, we may work with the smooth eigenvalue branch $h_{11}$ near $X_0$ and set
\[
   \widehat Q=\log h_{11}-\log(u-a)+\frac A2|X|^2.
\]
It has a local maximum at $X_0$.  From now on all computations are at $X_0$.

Since $G$ is symmetric, elliptic, and concave, the ordering above gives
\[
  G^{11}\ge G^{22}\ge\cdots\ge G^{nn},
  \qquad
  F^{11}\le F^{22}\le\cdots\le F^{nn}.
\]
The first- and second-derivative conditions are
\begin{equation}\label{eq:first-Q}
   0=\widehat Q_i=\frac{h_{11i}}{h_{11}}-\frac{u_i}{u-a}+A\langle X,e_i\rangle,
\end{equation}
and
\begin{equation}\label{eq:second-Q}
  0\ge F^{ii}\widehat Q_{ii}
  =F^{ii}(\log h_{11})_{ii}-F^{ii}(\log(u-a))_{ii}
   +\frac A2F^{ii}(|X|^2)_{ii}.
\end{equation}

\begin{lemma}\label{lem:master}
At $X_0$, if $h_{11}\ge1$, then
\begin{align}\label{eq:master}
0\ge{}&-\frac2{h_{11}}\sum_{i\ge2}G^{1i,i1}h_{11i}^2
       -\frac{F^{ii}h_{11i}^2}{h_{11}^2}
       +\frac{a}{u-a}F^{ii}h_{ii}^2
       +\frac{F^{ii}u_i^2}{(u-a)^2}\notag\\
&\quad -(C+R_\al)h_{11}+A\sum_iF^{ii}-CA.
\end{align}
\end{lemma}

\begin{proof}
From \eqref{eq:gauss}, \eqref{eq:trace-Fh},
\begin{align}\label{eq:X2-term}
 \frac A2F^{ii}(|X|^2)_{ii}
 &=A\sum_iF^{ii}(1-h_{ii}u)\notag\\
 &=A\sum_iF^{ii}-Au\ft+AuR_\al.
\end{align}
Moreover, by \eqref{eq:gauss}--\eqref{eq:codazzi},
\begin{equation}\label{eq:u-derivatives}
   u_i=h_{ii}\langle X,e_i\rangle,
   \qquad
   u_{ii}=\sum_kh_{iik}\langle X,e_k\rangle-uh_{ii}^2+h_{ii}.
\end{equation}
Thus
\begin{align}\label{eq:log-u}
 -F^{ii}(\log(u-a))_{ii}
 ={}&-\frac1{u-a}\sum_kF^{ii}h_{iik}\langle X,e_k\rangle
 +\frac{u}{u-a}F^{ii}h_{ii}^2\notag\\
 &-\frac{\ft-R_\al}{u-a}
 +\frac{F^{ii}u_i^2}{(u-a)^2}.
\end{align}
Since $\eta_{ii,k}=H_k-h_{ii,k}$ and \eqref{eq:Fii} holds,
\[
   F^{ii}h_{ii,k}=G^{ii}\eta_{ii,k}.
\]
Differentiating $G(\eta)=\ft$ in the direction $e_k$ gives
\[
   G^{ii}\eta_{ii,k}=\ft_k.
\]
At the diagonal point $X_0$, $\nu_k=h_{kk}e_k$, and the chain rule for
$\ft=\ft(X,\nu)$ therefore yields
\[
   \ft_k=(d_X\ft)(e_k)+h_{kk}(d_\nu\ft)(e_k).
\]
On the other hand,
\[
 \sum_iF^{ii}h_{ii,k}
 =\sum_i\sum_{p\ne i}G^{pp}h_{ii,k}
 =\sum_pG^{pp}(H_k-h_{pp,k})
 =\sum_pG^{pp}\eta_{pp,k}.
\]
Hence
\begin{equation}\label{eq:first-equation}
  F^{ii}h_{ii,k}
  =(d_X\ft)(e_k)+h_{kk}(d_\nu\ft)(e_k).
\end{equation}
Indeed, by \eqref{eq:first-equation},
\[
\begin{aligned}
	-\frac{1}{u-a}\sum_k F^{ii} h_{ii,k} \langle X,e_k\rangle
	&= -\frac{1}{u-a}\sum_k h_{kk} (d_\nu \widetilde{f})(e_k)\langle X,e_k\rangle \\
	&\quad -\frac{1}{u-a}(d_X \widetilde{f})(X^T),
\end{aligned}
\]
where $X^T=\sum_k \langle X,e_k\rangle e_k$. Moreover, $|X^T|\le |X|\le C$, $u-a\ge a>0$, and $|d_X\widetilde{f}|\le C$, because $d_X\widetilde{f}=\frac1k f^{1/k-1}d_X f$, $\inf f>0$, and $\|f\|_{C^1}\le C$. Hence
\[
-\frac{1}{u-a}\sum_k F^{ii} h_{ii,k} \langle X,e_k\rangle
\ge -\frac{1}{u-a}\sum_k h_{kk}(d_\nu \widetilde{f})(e_k)\langle X,e_k\rangle -C.
\]

Using \eqref{eq:trace-Fh}, we therefore obtain
\begin{align}\label{eq:log-u-lower}
 -F^{ii}(\log(u-a))_{ii}
 \ge{}&-\frac1{u-a}\sum_kh_{kk}(d_\nu\ft)(e_k)\langle X,e_k\rangle
 +\frac{u}{u-a}F^{ii}h_{ii}^2\notag\\
 &+\frac{R_\al}{u-a}+\frac{F^{ii}u_i^2}{(u-a)^2}-C.
\end{align}

For the first term in \eqref{eq:second-Q},
\begin{equation}\label{eq:log-h}
  F^{ii}(\log h_{11})_{ii}
  =\frac{F^{ii}h_{11ii}}{h_{11}}-\frac{F^{ii}h_{11i}^2}{h_{11}^2}.
\end{equation}
The commutation identity \eqref{eq:commutation-h} and \eqref{eq:trace-Fh} yield
\begin{equation}\label{eq:commuted}
  F^{ii}h_{11ii}
  =F^{ii}h_{ii11}-F^{ii}h_{ii}^2h_{11}+(\ft-R_\al)h_{11}^2.
\end{equation}
Differentiating $G(\eta)=\ft$ twice in the $e_1$ direction gives
\[
  F^{ii}h_{ii11}
  =-G^{ij,kl}\eta_{ij,1}\eta_{kl,1}+\ft_{11}.
\]
We next make the estimate for $\ft_{11}$ explicit.  Extend $f$ smoothly in the normal variable to a
neighborhood of $\mathbb S^n$; all derivatives of this extension are controlled by the given
$C^2$ norm.  At $X_0$ the frame is geodesic and $h$ is diagonal.  Applying the chain rule twice to
$\ft(X,\nu)$ produces a term
\[
   \sum_k h_{k11}(d_\nu\ft)(e_k)
\]
from the tangential part of $\nu_{11}$.  Every remaining term contains either a bounded first or
second derivative of $\ft$, or at most two factors of $h_{1j}$.  Since $h_{1j}=0$ for $j\ne1$ at
$X_0$, those remaining terms are bounded below by $-C(1+h_{11}^2)$.  Thus
\begin{equation}\label{eq:ft-second}
  F^{ii}h_{ii11}
  \ge -G^{ij,kl}\eta_{ij,1}\eta_{kl,1}
      +\sum_kh_{k11}(d_\nu\ft)(e_k)-Ch_{11}^2-C.
\end{equation}

For completeness, at a diagonal matrix the standard second-derivative formula for a symmetric
spectral function gives
\[
 G^{ij,kl}\xi_{ij}\xi_{kl}
 =\sum_{p,q}G^{pp,qq}\xi_{pp}\xi_{qq}
   +2\sum_{p<q}\frac{G^{pp}-G^{qq}}{\eta_{pp}-\eta_{qq}}\xi_{pq}^2,
\]
with the quotient understood by continuity when two eigenvalues coincide.  Concavity makes the
whole quadratic form nonpositive and, since $\eta_{11}\le\eta_{ii}$ and
$G^{11}\ge G^{ii}$, one has $G^{1i,i1}\le0$.  Keeping only the off-diagonal terms involving
index $1$, and using
\[
   \eta_{1i,1}=-h_{1i,1}=-h_{11i}
\]
by Codazzi, we obtain
\begin{equation}\label{eq:concavity-third}
 -G^{ij,kl}\eta_{ij,1}\eta_{kl,1}
 \ge -2\sum_{i\ge2}G^{1i,i1}h_{11i}^2.
\end{equation}
Substituting \eqref{eq:ft-second}--\eqref{eq:concavity-third} into
\eqref{eq:commuted} and then \eqref{eq:log-h} gives
\begin{align}\label{eq:log-h-lower}
F^{ii}(\log h_{11})_{ii}
\ge{}&-\frac2{h_{11}}\sum_{i\ge2}G^{1i,i1}h_{11i}^2
 -\frac{F^{ii}h_{11i}^2}{h_{11}^2}-F^{ii}h_{ii}^2\notag\\
&+\frac1{h_{11}}\sum_kh_{k11}(d_\nu\ft)(e_k)
 -(C+R_\al)h_{11}.
\end{align}
The two $d_\nu\ft$ terms in \eqref{eq:log-u-lower} and \eqref{eq:log-h-lower} actually
cancel up to the harmless $A$-term.  Indeed, Codazzi gives $h_{k11}=h_{11k}$, while
\eqref{eq:first-Q} and $u_k=h_{kk}\langle X,e_k\rangle$ give
\[
 \frac{h_{k11}}{h_{11}}
 =\frac{h_{kk}\langle X,e_k\rangle}{u-a}-A\langle X,e_k\rangle.
\]
Hence
\begin{align*}
&\frac1{h_{11}}\sum_kh_{k11}(d_\nu\ft)(e_k)
 -\frac1{u-a}\sum_kh_{kk}(d_\nu\ft)(e_k)\langle X,e_k\rangle\\
&\qquad=-A\sum_k\langle X,e_k\rangle(d_\nu\ft)(e_k)\ge-CA,
\end{align*}
where the last inequality uses the $C^1$ bound for $f$ and the $C^0$ bound for $|X|$.
Combining \eqref{eq:X2-term}, \eqref{eq:log-u-lower}, and
\eqref{eq:log-h-lower} in \eqref{eq:second-Q}, and discarding the nonnegative terms
$R_\al/(u-a)$ and $AuR_\al$, proves \eqref{eq:master}.
\end{proof}

We next absorb the unfavorable term in \eqref{eq:master}.

\begin{lemma}\label{lem:badterm}
If $h_{11}$ and $A$ are sufficiently large, depending only on the data in
Theorem~\ref{thm:curvature}, then
\begin{equation}\label{eq:badterm}
  (C+R_\al)h_{11}
  \le \frac{a}{2(u-a)}\sum_iF^{ii}h_{ii}^2
      +\frac A2\sum_iF^{ii}.
\end{equation}
\end{lemma}

\begin{proof}
Fix a small number $\delta>0$.

\emph{Case 1.} Suppose $|h_{ii}|\le\delta h_{11}$ for every $i\ge2$. Then
\[
   |\eta_{11}|\le(n-1)\delta h_{11}
\]
and
\begin{equation}\label{eq:eta-case1}
  [1-(n-2)\delta]h_{11}
  \le\eta_{22}\le\cdots\le\eta_{nn}
  \le[1+(n-2)\delta]h_{11}.
\end{equation}
Introduce the normalized eigenvalues $\mu_i=\eta_{ii}/h_{11}$.  By \eqref{eq:eta-case1}, under the present
assumption,
\[
 |\mu_1|\le(n-1)\delta,
 \qquad
 1-(n-2)\delta\le\mu_i\le1+(n-2)\delta\quad(i\ge2).
\]
Thus $\mu$ lies in a fixed small neighborhood of $(0,1,\ldots,1)$.  Since
\[
 \sig_{k-1}(0,1,\ldots,1)=\binom{n-1}{k-1}>0,
 \qquad
 \sig_{k-2}(0,1,\ldots,1)=\binom{n-1}{k-2}>0,
\]
continuity gives, after fixing $\delta>0$ sufficiently small,
\[
  \sig_{k-1}(\eta)=h_{11}^{k-1}\sig_{k-1}(\mu)\ge c h_{11}^{k-1},
  \qquad
  \sig_{k-2}(\eta)=h_{11}^{k-2}\sig_{k-2}(\mu)\ge c h_{11}^{k-2}.
\]
Using the exact trace formula from Lemma~\ref{lem:trace-lower},
\[
 \sum_iF^{ii}
 =\frac{n-1}{k}f^{1/k-1}
 \bigl[(n-k+1)\sig_{k-1}(\eta)
       +\al(n-k+2)\sig_{k-2}(\eta)\bigr],
\]
and the positive upper and lower bounds for $f$, we obtain
\[
   \sum_iF^{ii}\ge c h_{11}^{k-1}.
\]
Since $k\ge2$ and $h_{11}\ge1$, $h_{11}^{k-1}\ge h_{11}$.  We may therefore fix $A$
large enough, depending only on the data, so that
\[
   Ch_{11}\le\frac A2\sum_iF^{ii}.
\]

It remains to estimate $R_\al h_{11}$.  In terms of the normalized vector $\mu$, as
$\delta\to0$ one has
\[
   \mu\to(0,1,\ldots,1).
\]
Because $k\le n-1$, for every $j\ge2$,
\[
   \sig_{k-1}((0,1,\ldots,1)|j)=\binom{n-2}{k-1}>0,
\]
and
\[
   \sig_{k-2}((0,1,\ldots,1)|j)=\binom{n-2}{k-2}>0.
\]
For $j\ge2$,
\[
 G^{jj}=\frac1k f^{1/k-1}
 \bigl[\sig_{k-1}(\eta|j)+\al\sig_{k-2}(\eta|j)\bigr].
\]
The two displayed limiting values above are strictly positive, so, after decreasing the already
fixed $\delta$ if necessary,
\[
 \sig_{k-1}(\eta|j)\ge c h_{11}^{k-1},
 \qquad
 \sig_{k-2}(\eta|j)\ge c h_{11}^{k-2},
 \qquad j\ge2.
\]
Hence
\[
   G^{jj}\ge c f^{1/k-1}
   \bigl(h_{11}^{k-1}+\al h_{11}^{k-2}\bigr),\qquad j\ge2.
\]
On the other hand, all components of $\mu$ are uniformly bounded, so the same formula for
$G^{ii}$ gives
\[
  \sum_iG^{ii}\le C f^{1/k-1}
  \bigl(h_{11}^{k-1}+\al h_{11}^{k-2}\bigr).
\]
Since $\sum_iF^{ii}=(n-1)\sum_iG^{ii}$, it follows that
\[
    F^{11}=\sum_{j\ne1}G^{jj}\ge c\sum_iF^{ii}.
\]
It follows that
\begin{equation}\label{eq:quadratic-case1}
  \sum_iF^{ii}h_{ii}^2\ge F^{11}h_{11}^2
  \ge c h_{11}^2\sum_iF^{ii}.
\end{equation}
By \eqref{eq:Ralpha-trace}, $R_\al\le C\sum_iF^{ii}$, while
\eqref{eq:quadratic-case1} gives
\[
 \sum_iF^{ii}\le\frac{C}{h_{11}^2}\sum_iF^{ii}h_{ii}^2.
\]
Consequently
\[
 R_\al h_{11}
 \le\frac{C}{h_{11}}\sum_iF^{ii}h_{ii}^2.
\]
Because $a>0$ and $u$ is uniformly bounded above and below by \eqref{eq:u-bounds}, the
coefficient $a/[4(u-a)]$ has a fixed positive lower bound.  Thus, for sufficiently large
$h_{11}$,
\[
   R_\al h_{11}
   \le \frac{a}{4(u-a)}\sum_iF^{ii}h_{ii}^2.
\]
This proves \eqref{eq:badterm} in Case 1.

\emph{Case 2.} Suppose $h_{22}>\delta h_{11}$ or $h_{nn}<-\delta h_{11}$. Since
$G^{11}\ge G^{ii}$ and \eqref{eq:Fii} holds, for every $i\ge2$,
\begin{equation}\label{eq:F-lower-each}
  F^{ii}\ge G^{11}\ge\frac1n\sum_jG^{jj}
  =\frac1{n(n-1)}\sum_jF^{jj}.
\end{equation}
Therefore
\begin{equation}\label{eq:quadratic-case2}
   \sum_iF^{ii}h_{ii}^2
   \ge c\delta^2h_{11}^2\sum_iF^{ii}.
\end{equation}
Indeed, Lemma~\ref{lem:trace-lower} gives $\sum_iF^{ii}\ge c>0$, so
\[
 Ch_{11}\le C h_{11}\sum_iF^{ii},
 \qquad
 R_\al h_{11}\le C h_{11}\sum_iF^{ii}
\]
by \eqref{eq:Ralpha-trace}.  In view of \eqref{eq:quadratic-case2}, the right-hand sides
are bounded by
\[
 \frac{C}{\delta^2h_{11}}\sum_iF^{ii}h_{ii}^2.
\]
For sufficiently large $h_{11}$ this is at most
$a[4(u-a)]^{-1}\sum_iF^{ii}h_{ii}^2$.  This proves \eqref{eq:badterm} in Case 2 and finishes
the proof.
\end{proof}

Combining Lemmas~\ref{lem:master} and \ref{lem:badterm},
\begin{align}\label{eq:post-badterm}
0\ge{}&-\frac2{h_{11}}\sum_{i\ge2}G^{1i,i1}h_{11i}^2
       -\frac{F^{ii}h_{11i}^2}{h_{11}^2}
       +\frac{a}{2(u-a)}F^{ii}h_{ii}^2\notag\\
&\quad +\frac{F^{ii}u_i^2}{(u-a)^2}
       +\frac A2\sum_iF^{ii}-CA.
\end{align}

\begin{lemma}\label{lem:partial-curvature}
If $h_{11}$ and $A$ are sufficiently large, then at $X_0$,
\[
    |h_{ii}|\le CA,\qquad i\ge2.
\]
\end{lemma}

\begin{proof}
From \eqref{eq:first-Q}, for any $0<\varepsilon<1$,
\[
  \frac{F^{ii}h_{11i}^2}{h_{11}^2}
  \le\frac{1+\varepsilon}{(u-a)^2}F^{ii}u_i^2
   +C\left(1+\frac1\varepsilon\right)A^2F^{ii}\langle X,e_i\rangle^2.
\]
Substitute this into \eqref{eq:post-badterm} and discard the nonnegative concavity term
$-2h_{11}^{-1}\sum_{i\ge2}G^{1i,i1}h_{11i}^2$.  Since
$u_i=h_{ii}\langle X,e_i\rangle$, the excess part of the $u_i^2$ term is bounded by
\[
 \frac{\varepsilon}{(u-a)^2}F^{ii}u_i^2
 \le C\varepsilon F^{ii}h_{ii}^2,
\]
where we used $|X|\le C$ and $u-a\ge a>0$.  Choose $\varepsilon>0$ so small that this term
is absorbed by one half of the positive quadratic term
$a[2(u-a)]^{-1}F^{ii}h_{ii}^2$.  The remaining $A$-term satisfies
\[
 C\left(1+\frac1\varepsilon\right)A^2F^{ii}\langle X,e_i\rangle^2
 \le CA^2\sum_iF^{ii}.
\]
We therefore obtain
\begin{equation}\label{eq:partial-mid}
   0\ge c\sum_iF^{ii}h_{ii}^2
       -CA^2\sum_iF^{ii}-CA.
\end{equation}
Starting from \eqref{eq:partial-mid}, and using \eqref{eq:trace-lower}, $\sum_iF^{ii}\ge c_0>0$; since $A\ge1$, the final term satisfies
$CA\le C A^2\sum_iF^{ii}$ and can be absorbed into the middle term after enlarging its
constant.  Thus
\[
 c\sum_iF^{ii}h_{ii}^2\le CA^2\sum_iF^{ii}.
\]
For each $i\ge2$, \eqref{eq:F-lower-each} gives
$F^{ii}\ge c\sum_jF^{jj}$.  Hence
\[
 c\Bigl(\sum_jF^{jj}\Bigr)\sum_{i\ge2}h_{ii}^2
 \le CA^2\sum_jF^{jj}.
\]
Cancelling the positive trace gives
\[
   \sum_{i\ge2}h_{ii}^2\le CA^2,
\]
and therefore $|h_{ii}|\le CA$ for every $i\ge2$.
\end{proof}

\begin{proof}[Proof of Theorem~\ref{thm:curvature}]
If $h_{11}$ is bounded by a fixed multiple of $A$, there is nothing to prove. Hence we may
assume it is sufficiently large. By Lemma~\ref{lem:partial-curvature}, after increasing this
threshold we may also assume
\begin{equation}\label{eq:small-transverse}
   |h_{ii}|\le\delta h_{11},\qquad i\ge2,
\end{equation}
where $\delta>0$ is fixed sufficiently small.

For the $i=1$ derivative, \eqref{eq:first-Q} gives
\[
 \frac{h_{111}}{h_{11}}=\frac{u_1}{u-a}-A\langle X,e_1\rangle.
\]
For any $\varepsilon>0$, Cauchy--Schwarz yields
\[
 \frac{F^{11}h_{111}^2}{h_{11}^2}
 \le (1+\varepsilon)\frac{F^{11}u_1^2}{(u-a)^2}
 +C_\varepsilon A^2F^{11}\langle X,e_1\rangle^2.
\]
Since $u_1=h_{11}\langle X,e_1\rangle$, the excess term
$\varepsilon F^{11}u_1^2/(u-a)^2$ is bounded by
$C\varepsilon F^{11}h_{11}^2$.  Choose $\varepsilon$ sufficiently small; then, once
$h_{11}\gg A$, both this excess and the $A^2F^{11}$ term are absorbed by a fixed fraction of
$a(u-a)^{-1}\sum_iF^{ii}h_{ii}^2$.  Thus
\begin{equation}\label{eq:h111}
  \frac{F^{11}h_{111}^2}{h_{11}^2}
  \le \frac{F^{11}u_1^2}{(u-a)^2}
      +\frac{a}{8(u-a)}F^{ii}h_{ii}^2.
\end{equation}

For $i\ge2$, the standard second-derivative formula for a symmetric spectral function gives
\begin{equation*}
  -G^{1i,i1}=\frac{G^{11}-G^{ii}}{\eta_{ii}-\eta_{11}}
  =\frac{F^{ii}-F^{11}}{h_{11}-h_{ii}}.
\end{equation*}
Set $d_i=h_{11}-h_{ii}$.  By \eqref{eq:small-transverse},
\[
 (1-\delta)h_{11}\le d_i\le(1+\delta)h_{11},\qquad i\ge2.
\]
The spectral second-derivative formula therefore gives
\[
 -\frac2{h_{11}}G^{1i,i1}h_{11i}^2
 =\frac{2(F^{ii}-F^{11})}{h_{11}d_i}h_{11i}^2
 \ge\frac{2}{1+\delta}\frac{(F^{ii}-F^{11})h_{11i}^2}{h_{11}^2}.
\]
Since
\[
 F^{ii}\le\frac{2}{1+\delta}(F^{ii}-F^{11})
      +\frac{2}{1+\delta}F^{11},
\]
we obtain, after summing over $i\ge2$,
\[
 \sum_{i\ge2}\frac{F^{ii}h_{11i}^2}{h_{11}^2}
 \le -\frac2{h_{11}}\sum_{i\ge2}G^{1i,i1}h_{11i}^2
 +C F^{11}\sum_{i\ge2}\frac{h_{11i}^2}{h_{11}^2}.
\]
For $i\ge2$, \eqref{eq:first-Q}, $u_i=h_{ii}\langle X,e_i\rangle$, and
Lemma~\ref{lem:partial-curvature} give
\[
 \left|\frac{h_{11i}}{h_{11}}\right|
 \le \frac{|h_{ii}|}{u-a}|\langle X,e_i\rangle|
      +A|\langle X,e_i\rangle|
 \le CA|\langle X,e_i\rangle|.
\]
Hence
\[
 F^{11}\sum_{i\ge2}\frac{h_{11i}^2}{h_{11}^2}\le CA^2F^{11}.
\]
If $h_{11}$ is sufficiently large compared with $A$, this is bounded by
$a[8(u-a)]^{-1}F^{11}h_{11}^2$.  We have thus proved
\begin{align}\label{eq:third-order-control}
 \sum_{i\ge2}\frac{F^{ii}h_{11i}^2}{h_{11}^2}
 &\le -\frac2{h_{11}}\sum_{i\ge2}G^{1i,i1}h_{11i}^2
      +\frac{a}{8(u-a)}F^{11}h_{11}^2.
\end{align}

Substituting \eqref{eq:h111} and \eqref{eq:third-order-control} into
\eqref{eq:post-badterm} yields
\begin{equation}\label{eq:final-main}
  0\ge \frac{a}{4(u-a)}\sum_iF^{ii}h_{ii}^2
      +\frac A2\sum_iF^{ii}-CA.
\end{equation}
Dropping the first nonnegative term in \eqref{eq:final-main} gives
\[
 \frac A2\sum_iF^{ii}\le CA.
\]
Since $A>0$ is fixed, $\sum_iF^{ii}\le C$.  Using
$\sum_iF^{ii}=(n-1)\sum_iG^{ii}$, we obtain
\begin{equation}\label{eq:Gtrace-upper}
    \sum_iG^{ii}\le C.
\end{equation}

It remains to obtain a positive lower bound for $F^{11}$. Introduce the augmented vector
\[
   \widehat\eta=(\eta_{11},\ldots,\eta_{nn},\al)\in\Ga_k^{n+1}.
\]
Let $\widehat G=\sig_k(\widehat\eta)^{1/k}$. For the first $n$ coordinates,
$\widehat G^{ii}=G^{ii}$, while
\[
  \widehat G^{n+1,n+1}
  =\frac1k f^{1/k-1}\sig_{k-1}(\eta)
  \le\frac1{n-k+1}\sum_{i=1}^nG^{ii}.
\]
Thus \eqref{eq:Gtrace-upper} implies
\begin{equation}\label{eq:hatGtrace-upper}
   \sum_{a=1}^{n+1}\widehat G^{aa}\le C.
\end{equation}
Indeed, writing $N=n+1$, the identity
\[
  \sig_k(\widehat\eta)^{1/k-1}\sig_{k-1}(\widehat\eta)
  =\frac{k}{N-k+1}\sum_{a=1}^{N}\widehat G^{aa}\le C
\]
verifies the hypothesis of part~(2) of Lemma~2.2 in
Hou--Ma--Wu~\cite{HouMaWu2010}. Applied to
$\widehat\eta\in\Ga_k^{N}$, that lemma yields a positive lower bound for the
$\widehat G$-derivative corresponding to the largest component of
$\widehat\eta$.

By Lemma~\ref{lem:partial-curvature},
\[
  \eta_{11}=\sum_{j\ge2}h_{jj}=O(A),
\]
whereas, for $i\ge2$,
\[
  \eta_{ii}=h_{11}+O(A).
\]
Hence, if $h_{11}$ is sufficiently large, a largest component of $\widehat\eta$ is some
$\eta_{pp}$ with $p\ge2$. Consequently,
\[
   G^{pp}\ge c>0,
\]
and therefore
\begin{equation}\label{eq:F11-lower-final}
   F^{11}=\sum_{p\ne1}G^{pp}\ge c>0.
\end{equation}
Using \eqref{eq:F11-lower-final} in the first term of \eqref{eq:final-main}, and using the
uniform bounds for $u-a$, gives
\[
   0\ge cF^{11}h_{11}^2-CA\ge c h_{11}^2-CA.
\]
Hence $h_{11}\le C$.  Finally, since $H>0$ and
$h_{11}\ge h_{22}\ge\cdots\ge h_{nn}$,
\[
 0<H=\sum_{i=1}^nh_{ii}\le(n-1)h_{11}+h_{nn},
\]
so $h_{nn}\ge-(n-1)h_{11}$.  Thus every principal curvature is bounded in absolute value,
and the theorem follows.
\end{proof}

\section{A priori estimates and proof of Theorem~\ref{thm:existence}}\label{sec:existence}

We prove Theorem~\ref{thm:existence} by the continuity method.  In order to make the
argument transparent, we first establish the $C^0$, $C^1$, and $C^2$ estimates separately,
and then treat existence and uniqueness in the last two subsections.

Set
\[
   A_{n,k}=\binom nk(n-1)^k,
   \qquad
   B_{n,k,\al}=\al\binom n{k-1}(n-1)^{k-1},
\]
so that
\[
   \Phi_\al(r)=\frac{A_{n,k}}{r^k}+\frac{B_{n,k,\al}}{r^{k-1}}.
\]
Let
\[
   C_{n,k,\al}=A_{n,k}+B_{n,k,\al}.
\]
Choose
\[
   0<\varepsilon<\frac1{r_2^k-1}
\]
and define
\[
   f_0(X,\nu)
   =C_{n,k,\al}\left(\frac{1+\varepsilon}{|X|^k}-\varepsilon\right).
\]
The choice of $\varepsilon$ gives
\[
 \frac{1+\varepsilon}{r_2^k}-\varepsilon
 =\frac{1-\varepsilon(r_2^k-1)}{r_2^k}>0.
\]
Since $\rho\mapsto(1+\varepsilon)\rho^{-k}-\varepsilon$ is decreasing, it follows that
$f_0>0$ throughout $r_1\le|X|\le r_2$.  On the unit sphere,
$\lambda_i(\eta)=n-1$ for every $i$, and therefore
\[
 S_{k,\al}(\lambda(\eta))
 =A_{n,k}+B_{n,k,\al}=C_{n,k,\al}=f_0.
\]
Thus the unit sphere is a solution of the $t=0$ equation.  Moreover,
\[
   \frac{\partial}{\partial\rho}(\rho^kf_0)
   =-k\varepsilon C_{n,k,\al}\rho^{k-1}<0.
\]
A direct calculation gives
\[
 \rho^k\bigl(f_0(\rho)-\Phi_\al(\rho)\bigr)
 =B_{n,k,\al}(1-\rho)
  +\varepsilon C_{n,k,\al}(1-\rho^k).
\]
Thus
\[
   f_0(\rho)>\Phi_\al(\rho)\quad(\rho<1),
   \qquad
   f_0(\rho)<\Phi_\al(\rho)\quad(\rho>1).
\]
For $t\in[0,1]$, define
\[
    f_t=tf+(1-t)f_0.
\]
Because the path is a convex combination, $f_t$ remains positive.  The quantity
$\partial_\rho(\rho^kf_t)$ is the same convex combination of
$\partial_\rho(\rho^kf)$ and $\partial_\rho(\rho^kf_0)$; hence
\eqref{eq:radial-monotonicity} holds for every $t$.  Likewise, $f$ satisfies the barrier inequalities in \eqref{eq:barrier}
non-strictly while $f_0$ satisfies them strictly at $r_1$ and $r_2$.
Therefore, for every $t<1$, the two barrier inequalities for $f_t$ are strict.  All constants
below are independent of $t$.

\subsection{\texorpdfstring{$C^0$ estimate}{C0 estimate}}\label{subsec:C0-estimate}

Let $M_t$ be an admissible star-shaped solution with radial function $\rho$.  At a maximum
point of $\rho$,
\[
  \nabla\rho=0,\qquad \nabla^2\rho\le0.
\]
By \eqref{eq:radial-formulas},
\[
   g_{ij}=\rho^2\widehat g_{ij},
   \qquad
   h_{ij}=\rho\widehat g_{ij}-\rho_{ij}\ge\rho\widehat g_{ij},
\]
and hence
\[
   \kappa_i\ge\frac1\rho,
   \qquad
   \lambda_i(\eta)\ge\frac{n-1}{\rho}.
\]
Since $S_{k,\al}$ is increasing in each eigenvalue throughout the admissible region, and the
admissible region is convex by Lemma~\ref{lem:augmented}, componentwise comparison along the line
segment joining $\lambda(\eta)$ to the spherical vector gives
\begin{equation}\label{eq:C0-max}
   f_t(X,\nu)=S_{k,\al}(\lambda(\eta))\ge\Phi_\al(\rho).
\end{equation}
At such an extremum $\nu=x=X/|X|$.  If a solution in the closure of the annular class touched
$\rho=r_2$ at its maximum, then \eqref{eq:C0-max} would give
$f_t(r_2x,x)\ge\Phi_\al(r_2)$, contradicting the strict outer barrier for $t<1$.

At a minimum point of $\rho$ the inequalities reverse:
\[
   h_{ij}\le\rho\widehat g_{ij},
   \qquad
   \lambda_i(\eta)\le\frac{n-1}{\rho},
\]
and therefore
\[
   f_t(X,\nu)\le\Phi_\al(\rho).
\]
Again $\nu=x$ at the extremum.  If $\rho=r_1$ at the minimum, this inequality contradicts the
strict inner barrier $f_t(r_1x,x)>\Phi_\al(r_1)$ for $t<1$.  Consequently,
\begin{equation}\label{eq:C0}
   r_1<\rho<r_2\qquad (0\le t<1),
\end{equation}
with the corresponding non-strict bounds at $t=1$ by passage to the limit.

\subsection{\texorpdfstring{$C^1$ estimate}{C1 estimate}}\label{subsec:C1-estimate}

By \eqref{eq:support-rho}, it suffices to prove a positive lower bound for the support
function $u=\langle X,\nu\rangle$.  Consider
\[
    w=-\log u+\gamma(|X|^2),
    \qquad
    \gamma(s)=\frac{L}{s},
\]
where $L>1$ will be chosen sufficiently large.  Suppose that $w$ attains its maximum at
$X_0$.  If $X_0$ is parallel to $\nu$, then $u=|X|\ge r_1$, and there is nothing to prove.
Otherwise $X^T\ne0$.  Choose
$e_1=X^T/|X^T|$ and complete it to an orthonormal frame.  Then
\[
  \langle X,e_1\rangle\ne0,
  \qquad
  \langle X,e_i\rangle=0\quad(i\ge2).
\]
At the maximum point,
\[
 0=w_i=-\frac{u_i}{u}+2\gamma'\langle X,e_i\rangle.
\]
Since $u_i=h_{ij}\langle X,e_j\rangle$ and only the $e_1$ component of $X^T$ is nonzero,
we have $u_i=h_{i1}\langle X,e_1\rangle$.  Thus the $i=1$ equation gives
$h_{11}=2\gamma'u$, whereas the equations with $i\ge2$ give $h_{1i}=0$.  We then rotate
$e_2,\ldots,e_n$ to diagonalize the remaining $(n-1)\times(n-1)$ block without changing
$e_1$.
Since $\gamma'<0$, we have $h_{11}<0$.  Because
$H>0$, there is an index $p\ge2$ such that $h_{pp}>0$.  Hence $\eta_{pp}=H-h_{pp}<H-h_{11}=\eta_{11}$.  For a symmetric concave function, the
derivatives are ordered oppositely to the eigenvalues; hence $G^{pp}\ge G^{11}$.  Since
$G^{pp}$ is one of the terms in $F^{11}=\sum_{j\ne1}G^{jj}$, we obtain
$F^{11}\ge G^{pp}\ge G^{11}$.  Consequently
$\sum_jG^{jj}=G^{11}+F^{11}\le2F^{11}$, and therefore
\begin{equation}\label{eq:F11-gradient}
   F^{11}=\sum_{j\ne1}G^{jj}
   \ge\frac12\sum_jG^{jj}
   =\frac1{2(n-1)}\sum_iF^{ii}.
\end{equation}

At $X_0$,
\[
 w_{ii}=-\frac{u_{ii}}u+\frac{u_i^2}{u^2}
        +\gamma''\bigl((|X|^2)_i\bigr)^2
        +\gamma'(|X|^2)_{ii}.
\]
Because $(|X|^2)_i=2\langle X,e_i\rangle$ and only the $i=1$ tangential component is
nonzero, the first derivative relation gives
\[
 \frac{u_1^2}{u^2}=4\gamma'^2\langle X,e_1\rangle^2,
 \qquad u_i=0\quad(i\ge2).
\]
Multiplying by $F^{ii}$ and using $F^{ii}w_{ii}\le0$ yields
\begin{equation}\label{eq:w-second}
0\ge -\frac{F^{ii}u_{ii}}u
 +4(\gamma'^2+\gamma'')F^{11}\langle X,e_1\rangle^2
 +\gamma'F^{ii}(|X|^2)_{ii}.
\end{equation}
By \eqref{eq:u-derivatives},
\[
 F^{ii}u_{ii}
 =\sum_k\langle X,e_k\rangle F^{ii}h_{ii,k}
   -uF^{ii}h_{ii}^2+\sum_iF^{ii}h_{ii}.
\]
Only $k=1$ contributes to the first sum.  Using \eqref{eq:first-equation} and
\eqref{eq:trace-Fh} therefore gives
\begin{align*}
 F^{ii}u_{ii}
 ={}&\langle X,e_1\rangle(d_X\ft)(e_1)
 +h_{11}\langle X,e_1\rangle(d_\nu\ft)(e_1)\\
 &-uF^{ii}h_{ii}^2+\ft-R_\al.
\end{align*}
Similarly, from $(|X|^2)_{ii}=2(1-uh_{ii})$,
\[
 F^{ii}(|X|^2)_{ii}
 =2\sum_iF^{ii}-2u(\ft-R_\al).
\]
The radial condition for $f_t$ is
$k f_t+(d_Xf_t)(X)\le0$.  Since $\ft=f_t^{1/k}$,
\[
 \ft+(d_X\ft)(X)
 =\frac1k f_t^{1/k-1}\bigl(k f_t+(d_Xf_t)(X)\bigr)\le0.
\]
Since at $X_0$,
\[
   X=\langle X,e_1\rangle e_1+u\nu,
\]
the preceding radial inequality gives
\[
 -\frac{\langle X,e_1\rangle(d_X\ft)(e_1)+\ft}{u}
 \ge(d_X\ft)(\nu).
\]
Also $h_{11}=2\gamma'u$, so the normal-derivative contribution satisfies
\[
 -\frac{h_{11}\langle X,e_1\rangle}{u}(d_\nu\ft)(e_1)
 =-2\gamma'\langle X,e_1\rangle(d_\nu\ft)(e_1),
\]
which is bounded below by $-C|\gamma'|$.  The terms
$F^{ii}h_{ii}^2$ and $R_\al/u$ in $-F^{ii}u_{ii}/u$ are nonnegative and may be discarded.
Finally,
\[
 \gamma'F^{ii}(|X|^2)_{ii}
 =2\gamma'\sum_iF^{ii}-2u\gamma'\ft+2u\gamma'R_\al.
\]
Since $\gamma'<0$, the middle term is nonnegative, whereas
\eqref{eq:Ralpha-trace} bounds the last term below by
$-C\al|\gamma'|\sum_iF^{ii}$.  The remaining first derivatives of $\ft$ are uniformly
bounded because $f_t$ has a uniform positive lower bound and a uniform $C^1$ norm.
Substituting all these estimates into \eqref{eq:w-second} yields
\begin{equation}\label{eq:w-final}
  0\ge 4(\gamma'^2+\gamma'')F^{11}\langle X,e_1\rangle^2
    -C|\gamma'|(1+\al)\sum_iF^{ii}-C(1+|\gamma'|).
\end{equation}
If
\[
   \langle X,e_1\rangle^2<\frac12r_1^2,
\]
then $u^2=|X|^2-\langle X,e_1\rangle^2\ge r_1^2/2$, and again there is nothing to prove.
Otherwise, since $s=|X|^2$ lies in the fixed interval $[r_1^2,r_2^2]$ and
$\gamma(s)=L/s$,
\[
 \gamma'=-\frac{L}{s^2},\qquad \gamma''=\frac{2L}{s^3},
\]
so
\[
   \gamma'^2+\gamma''\ge cL^2,
   \qquad
   |\gamma'|\le CL.
\]
Using \eqref{eq:F11-gradient} and Lemma~\ref{lem:trace-lower}, \eqref{eq:w-final} becomes
\[
  0\ge [cL^2-CL(1+\al)]\sum_iF^{ii}-C(1+L).
\]
Because $\sum_iF^{ii}\ge c_{n,k}>0$, the right-hand side is positive for $L$ sufficiently
large, a contradiction.  Thus $u\ge c>0$.  Finally, from
\[
 u=\frac{\rho^2}{\sqrt{\rho^2+|\nabla\rho|^2}}
\]
and the $C^0$ bounds on $\rho$,
\[
 |\nabla\rho|^2=\frac{\rho^4}{u^2}-\rho^2\le C.
\]
Hence
\begin{equation}\label{eq:C1}
   |\nabla\rho|\le C.
\end{equation}

\subsection{\texorpdfstring{$C^2$ estimate}{C2 estimate}}\label{subsec:C2-estimate}

The estimates \eqref{eq:C0} and \eqref{eq:C1} give uniform upper and lower bounds for $\rho$
and a uniform bound for $|\nabla\rho|$.  By \eqref{eq:radial-formulas} this controls the first
fundamental form, its inverse, the support function, and the $C^1$ geometry entering the constant
in Theorem~\ref{thm:curvature}.  Moreover, since $f_t=tf+(1-t)f_0$ and the annulus
$r_1\le|X|\le r_2$ is fixed,
\[
   \inf f_t\ge c>0,
   \qquad
   \|f_t\|_{C^2}\le C
\]
uniformly for $t\in[0,1]$.  Hence all quantities on which the curvature estimate depends are
uniform along the path.  Theorem~\ref{thm:curvature} therefore yields
\[
   \max_{M_t}\max_{1\le i\le n}|\kappa_i|\le C.
\]
Thus the second fundamental form of $M_t$ is uniformly bounded.

It remains to translate this geometric curvature bound into a $C^2$ estimate for the radial
function.  Writing
\[
   W=\sqrt{\rho^2+|\nabla\rho|^2},
\]
the radial graph formula \eqref{eq:radial-formulas} gives
\[
   \rho_{ij}
   =\rho\widehat g_{ij}
     +\frac{2}{\rho}\rho_i\rho_j
     -\frac{W}{\rho}h_{ij}.
\]
The $C^0$ and $C^1$ estimates make $g$ and $\widehat g$ uniformly equivalent, so a bound for the
second fundamental form with respect to $g$ gives a bound for its coordinate components in a
$\widehat g$-orthonormal frame.  The factors $W/\rho$ and $2/\rho$ in the preceding identity are
also uniformly bounded.  Hence
\[
   |\nabla^2\rho|\le C,
\]
and therefore
\[
   \|\rho\|_{C^2(\Sn)}\le C
\]
uniformly for all admissible solutions along the path.

We also record why this gives uniform ellipticity.  The curvature bound makes
$\widehat\lambda=(\lambda(\eta),\al)$ uniformly bounded.  Moreover
$\sig_k(\widehat\lambda)=f_t\ge c>0$.  Since $\widehat\lambda\in\Gamma_k^{n+1}$,
Newton--Maclaurin gives positive lower bounds for the lower symmetric functions
$\sig_j(\widehat\lambda)$, $1\le j<k$, in terms of the lower bound for $\sig_k$.  Thus
$\widehat\lambda$ stays a positive distance away from $\partial\Gamma_k^{n+1}$ and hence in
a compact subset of $\Gamma_k^{n+1}$.  The linearized coefficients are therefore bounded above
and below uniformly, so the equation is uniformly elliptic along the continuity path.

\subsection{Existence}\label{subsec:existence}

Let
\[
  \mathcal I=\{t\in[0,1]: S_{k,\al}(\lambda(\eta))=f_t(X,\nu)
  \text{ has an admissible star-shaped solution}\}.
\]
By construction, $0\in\mathcal I$, with the unit sphere as the corresponding solution.

We first prove openness.  Write $\rho=e^v$.  The principal curvatures have the form
\[
   \kappa_i=e^{-v}\overline\kappa_i(Dv,D^2v),
\]
and therefore
\[
   \lambda_i(\eta)=e^{-v}\overline\lambda_i(Dv,D^2v).
\]
After multiplying the equation by $e^{kv}$, we obtain
\[
   \sig_k(\overline\lambda)
   +\al e^v\sig_{k-1}(\overline\lambda)
   =e^{kv}f_t(e^vx,\nu(Dv)).
\]
Define
\[
 \mathcal P[v]
 :=e^{kv}f_t(e^vx,\nu(Dv))
   -\sig_k(\overline\lambda)
   -\al e^v\sig_{k-1}(\overline\lambda).
\]
For fixed $(v,Dv)$, the radial graph formula shows that the matrix determining the principal
curvatures depends affinely on $-D^2v$ through a positive definite change of basis.  Since
$S_{k,\al}^{1/k}$ is elliptic on the admissible set, the above sign convention makes
$\mathcal P$ elliptic with respect to $D^2v$.  Holding $Dv$ and $D^2v$ fixed, we have
\begin{align}\label{eq:Pv}
 \mathcal P_v
 &=e^{kv}\bigl(kf_t+(d_Xf_t)(X)\bigr)
   -\al e^v\sig_{k-1}(\overline\lambda)\notag\\
 &=\rho\frac{\partial}{\partial\rho}(\rho^kf_t)
   -\al e^v\sig_{k-1}(\overline\lambda)<0.
\end{align}
The strict inequality follows from \eqref{eq:radial-monotonicity}, $\al>0$, and admissibility:
positive scaling preserves $\Gamma_{k-1}$, so
$\sig_{k-1}(\overline\lambda)>0$.  At a fixed solution the linearization therefore has the form
\[
 L\varphi=a^{ij}\nabla_{ij}\varphi+b^i\nabla_i\varphi+c\varphi,
 \qquad (a^{ij})>0,\qquad c<0.
\]
The coefficients are uniformly elliptic at that solution because admissibility is an open
condition.  If $L\varphi=0$, the maximum principle applied to the positive maximum and negative
minimum of $\varphi$ gives $\varphi\equiv0$; hence the kernel is trivial.  As a second-order
elliptic operator on the closed manifold $\mathbb S^n$, $L:C^{2,\beta}\to C^{0,\beta}$ is
Fredholm of index zero, so trivial kernel implies surjectivity and therefore invertibility.
The implicit function theorem then shows that $\mathcal I$ is open.

To prove closedness, let $t_j\in\mathcal I$ and $t_j\to t_\infty$, and write the corresponding
solutions as $\rho_j=e^{v_j}$.  The estimates in
Subsections~\ref{subsec:C0-estimate}--\ref{subsec:C2-estimate} give a uniform $C^2$ bound for
$v_j$ and keep the augmented eigenvalues in a fixed compact subset of $\Gamma_k^{n+1}$.
Thus the equations are uniformly elliptic with uniform ellipticity constants.  Rewrite the equation in the equivalent form
\[
  \Bigl[\sig_k(\overline\lambda)+\al e^v\sig_{k-1}(\overline\lambda)\Bigr]^{1/k}
  =e^v f_t(e^vx,\nu(Dv))^{1/k}.
\]
For fixed $(x,v,Dv)$, the left-hand side is concave in $D^2v$: it is the
restriction of $\sig_k^{1/k}$ to the augmented vector
$(\overline\lambda,\al e^v)$, with $\al e^v$ fixed, and the curvature matrix
is affine in $-D^2v$ after a positive change of basis.  The left-hand side is
decreasing with respect to $D^2v$ in the sense of the ordering of symmetric
matrices.  Thus, with the standard convention that elliptic operators are
increasing in $D^2v$, we use the equivalent equation $\mathcal Q[v]=0$, where
\[
  \mathcal Q[v]
  :=e^v f_t(e^vx,\nu(Dv))^{1/k}
   -\Bigl[\sig_k(\overline\lambda)
       +\al e^v\sig_{k-1}(\overline\lambda)\Bigr]^{1/k}.
\]
The operator $\mathcal Q$ is convex and uniformly elliptic in $D^2v$ on the
uniformly admissible solution family.  The Evans--Krylov regularity theory
\cite{Evans1982,Krylov1983} therefore gives
\[
 \|v_j\|_{C^{2,\theta}(\mathbb S^n)}\le C
\]
for some $0<\theta<1$.  Differentiating the equation gives a uniformly elliptic linear equation
for the first derivatives of $v_j$ whose coefficients are $C^\theta$; Schauder estimates then
give the $C^{3,\theta}$ bound needed for compactness (and, by the assumed $C^2$ regularity of
$f_t$, the stated $C^{3,\beta}$ regularity for every $\beta<1$ after the standard bootstrap on
slightly smaller H\"older exponents).  Hence a subsequence converges in $C^{3,\beta'}$ for every
$\beta'<\beta$ to a solution at $t=t_\infty$.  The compact containment in the admissible cone
shows that admissibility is preserved in the limit.  Therefore $\mathcal I$ is closed, and hence
\[
   \mathcal I=[0,1],
\]
and in particular a solution exists for $t=1$.

\subsection{Uniqueness}\label{subsec:uniqueness}

Let $v_1$ and $v_2$ be two admissible solutions of the equation at $t=1$.  Suppose, toward a
contradiction, that $w=v_1-v_2$ has a positive maximum at $x_0$.  Then
\[
   Dv_1(x_0)=Dv_2(x_0),
   \qquad
   D^2v_1(x_0)\le D^2v_2(x_0)
\]
as quadratic forms.  Keeping the common gradient and the value $v_1(x_0)$ fixed, ellipticity of
$\mathcal P$ in the Hessian variable gives
\[
  \mathcal P(v_1,Dv_1,D^2v_1)
  \le \mathcal P(v_1,Dv_2,D^2v_2).
\]
Now keep $(Dv_2,D^2v_2)$ fixed.  The barred eigenvalues
$\overline\lambda(Dv_2,D^2v_2)$ are independent of the zeroth-order variable $v$ and belong to
$\Gamma_{k-1}$ because $v_2$ is admissible.  Hence \eqref{eq:Pv} is strictly negative for every
$v$ between $v_2(x_0)$ and $v_1(x_0)$.  Since $v_1(x_0)>v_2(x_0)$,
\[
  \mathcal P(v_1,Dv_2,D^2v_2)
  <\mathcal P(v_2,Dv_2,D^2v_2).
\]
Both endpoints are solutions, so altogether
\[
  0=\mathcal P[v_1](x_0)
  <\mathcal P[v_2](x_0)=0,
\]
a contradiction.  Therefore $v_1\le v_2$.  Reversing the roles of the two solutions gives
$v_2\le v_1$, and hence $v_1=v_2$.  This proves uniqueness and completes the proof of
Theorem~\ref{thm:existence}.

\medskip
\noindent\textbf{Acknowledgments.}
Peihe Wang is supported by National Natural Science Foundation of Shandong Province(ZR2026MS0039). Hua Zhu is supported by the National Natural Science Foundation of China
(Grant Nos.~12501273 and 12661041).

\end{document}